\documentclass{amsart}

\usepackage{tristan}
\usepackage{calc}
\usepackage{pdfsync}
\usepackage{amsmath,amssymb}
\usepackage{enumerate}
\usepackage{xspace}
\usepackage{boxedminipage}
\usepackage{algorithmic}
\usepackage{listings}
\usepackage{tabularx}
\usepackage{subfigure}
\usepackage[english]{babel}
\usepackage{tikz}
\usepackage{fullpage}

\renewcommand{\bih}[2]{\ensuremath{\cA_h\qp{#1,#2}}}

\renewcommand{\R}{E^{\cR}}
\newcommand{\HCT}{\text{HCT}(k+2)}
\newcommand{\tv}{\widetilde{v}}
\newcommand{\tw}{\widetilde{w}}
\newcommand{\tg}{\widetilde{g}}

\newcommand{\dx}{\,{\rm d}x}
\newcommand{\dS}{\,{\rm d}s}
\renewcommand{\Norm}[1]{\|{#1}\|}

\renewcommand{\eenorm}[1]{\| {{#1}} \|_{2,h}}
\renewcommand{\enorm}[1]{\| {{#1}} \|_{1,h}}
\newcommand{\znorm}[1]{\| {{#1}} \|_{0,h}}

\begin{document}
\title[Analysis of DG methods using mesh dependent norms]{Analysis of Discontinuous Galerkin Methods \\ using Mesh-Dependent Norms \\ and Applications to Problems with Rough Data}

\author[EG]{Emmanuil H.~Georgoulis}

\author[TP]{Tristan Pryer}

\address[EG]{
  Department of Mathematics,
  University of Leicester,
  University Road,
  Leicester LE1 7RH,
  UK and Department of Mathematics, School of Applied Mathematical and Physical Sciences, National Technical University of Athens, Zografou 157 80, Greece.
  {\tt Emmanuil.Georgoulis@le.ac.uk}
}
\address[TP]{
  Department of Mathematics and Statistics,
  University of Reading,
  Whiteknights,
  PO Box 220,
  Reading RG6 6AX,
  UK.
  {\tt t.pryer@reading.ac.uk}
}


\begin{abstract}
We prove the inf-sup stability of a discontinuous Galerkin scheme for second order elliptic operators in (unbalanced) mesh-dependent norms for quasi-uniform meshes for all spatial dimensions. This results in a priori error bounds in these norms. As an application we examine a problem with rough source term where the solution can not be characterised as a weak solution and show quasi-optimal error control.
\end{abstract}

\maketitle

\section{Introduction}

Discontinuous Galerkin (dG) methods are a popular family of non-conforming finite element-type
approximation schemes for partial differential equations (PDEs) involving discontinuous approximation spaces. In the context of elliptic problems their inception can be traced back to the 1970s \cite{Nitsche:1971,Baker:1977}; see also \cite{ArnoldBrezziCockburnMarini:2001} for an accessible
overview and history of these methods for second order problems. For
higher order problems, for example the (nonlinear) biharmonic problem,
dG methods are a useful alternative to using $\cont{1}$-conforming
elements whose derivation and implementation can become very
complicated
\cite{Baker:1977,GeorgoulisHoustonVirtanen:2011,Pryer:2014}. 

Inf-sup conditions form one part of the Banach--Ne\v{c}as--Babu\v{s}ka condition which guarantees the well-posedness of a given variational problem. In this note, we shall describe an analytical framework to examine the stability of dG approximations for $\leb{2}$ and $\sobh{2}$-like mesh-dependent norms. This is in keeping with the spirit of \cite{BabuskaOsborn:1980,BabuskaOsbornPitkaranta:1980}, where for continuous finite element methods the authors prove equivalent results for second and fourth order problems respectively. The present approach, however, is quite different and results in inf-sup stability for both $\leb{2}$- and $\sobh{2}$-like mesh-dependent norms under the assumption that the underlying mesh is quasi-uniform.

The analysis presented utilises a new \emph{$\sobh{2}$-conforming reconstruction operator}, based on Hsieh--Clough--Tocher-type $\cont{1}$ reconstructions. Such reconstructions, based on nodal averaging, are used for the proof of a posteriori bounds for non-conforming methods for elliptic \cite{KarakashianPascal:2003,brenner_gudi_sung,GeorgoulisHoustonVirtanen:2011,Pryer:2015} and hyperbolic problems \cite{GeorgoulisHallMakridakis:2014,GiesselmannMakridakisPryer:2015}. The new reconstruction operators presented below enjoys certain orthogonality properties; in particular, they are \emph{adjoint orthogonal} to the underlying Hsieh--Clough--Tocher space and maintain the same stability bounds as the $\sobh{2}$-conforming reconstruction from \cite{GeorgoulisHoustonVirtanen:2011}.

The argument is quite general and allows the derivation of inf-sup stability results whenever the numerical scheme has a well posed discrete adjoint (dual) problem over an appropriately constructed  non-conforming finite element space. This is contrary to the Aubin--Nitsche $\leb{2}$ duality argument whereby it is the underlying partial differential operator itself that requires the well posedness of the adjoint \emph{continuous} problem.

The use of these recovery operators is not limited to an a posteriori setting, indeed, they have been used to quantify inconsistencies appearing in standard interior penalty methods when the exact solution is not $\sobh{2}(\W)$ \cite{Gudi:2010}. This allows for quasi optimal a priori bounds for elliptic problems under minimal regularity up to data oscillation. Fundamentally the assumption in this analysis is that the singularity arises from the geometry of the domain rather than through the problem data itself. Our analysis allows us to show quasi-optimal $\leb{2}$ convergence to problems that have rough problem data. To showcase the result we study the convergence of a method posed for an elliptic problem whose source term is not $\sobh{-1}$. In this case the Aubin--Nitsche and indeed the standard treatment of Galerkin methods are not applicable.

The note is set out as follows: In \S\ref{sec:setup} we introduce the problem and present the analysis cumulating in inf-sup stability for problems with smooth data. In \S\ref{sec:roughdata} we examine a particular problem with rough data and prove quasi-optimal convergence in this case. In addition we give some numerical validation of the method.

\section{Problem set up and discretisation}
\label{sec:setup}

To highlight the main steps of the present developments in this area, we consider the Poisson problem with homogeneous Dirichlet boundary conditions. Let $\W \subset \reals^d$ be a convex domain and consider the problem: Given $f\in\leb{2}(\W)$ find $u\in\sobh{2}(\W)\cap\hoz(\W)$, such that
\begin{equation}
  \label{eq:bilinear-form}
  \int_\W \nabla u \cdot \nabla v \dx
  =
  \int_\W f \ v \dx \Foreach v\in\hoz(\W).
\end{equation}
We consider $\T{}$ to be a conforming triangulation of $\W$,
namely, $\T{}$ is a finite family of sets such that
\begin{enumerate}
\item $K\in\T{}$ implies $K$ is an open simplex (segment for $d=1$,
  triangle for $d=2$, tetrahedron for $d=3$),
\item for any $K,J\in\T{}$ we have that $\closure K\meet\closure J$ is
  a full lower-dimensional simplex (i.e., it is either $\emptyset$, a vertex, an
  edge, a face, or the whole of $\closure K$ and $\closure J$) of both
  $\closure K$ and $\closure J$ and
\item $\union{K\in\T{}}\closure K=\closure\W$.
\end{enumerate}
The shape regularity constant of $\T{}$ is defined as the number
\begin{equation}
  \label{eqn:def:shape-regularity}
  \mu(\T{}) := \inf_{K\in\T{}} \frac{\rho_K}{h_K},
\end{equation}
where $\rho_K$ is the radius of the largest ball contained inside
$K$ and $h_K$ is the diameter of $K$. An indexed family of
triangulations $\setof{\T n}_n$ is called \emph{shape regular} if 
\begin{equation}
  \label{eqn:def:family-shape-regularity}
  \mu:=\inf_n\mu(\T n)>0.
\end{equation}
Further, we define $\funk h\W\reals$ to be the {piecewise
  constant} \emph{meshsize function} of $\T{}$ given by
\[
  h(\vec{x}):=\max_{\closure K\ni \vec{x}}h_K.
\]
A mesh is called quasiuniform when there exists a positive constant $C$ such that $\max_{x\in\Omega} h \le C \min_{x\in\Omega} h$. In what follows we shall assume that all triangulations are shape-regular and quasiuniform.

We let $\E{}$ be the skeleton (set of common interfaces) of the
triangulation $\T{}$ and say $e\in\E$ if $e$ is on the interior of
$\W$ and $e\in\partial\W$ if $e$ lies on the boundary $\partial\W$
{and set $h_e$ to be the diameter of $e$.} We also define the ``broken'' gradient $\nabla_h$, Laplacian $\Delta_h$ and Hessian $\Hess_h$ to be defined element-wise by $\nabla_h w|_K = \nabla w$, $\Delta_h w|_K = \Delta w$, $\Hess_h w|_K = \Hess w$ for all $K\in \T{}$, respectively, for respectively smooth functions on the interior of $K$, 

We let $\poly k(\T{})$ denote the space of piecewise polynomials of
degree $k$ over the triangulation $\T{}$,
 and introduce the \emph{finite element space}
 $\fes := \poly k(\T{})$
to be the usual space of discontinuous piecewise polynomial
functions of degree $k$. We define average operators for arbitrary scalar functions $v$ and vectors $\vec v$  over an edge $e$ shared by elements $K_1$ and $K_2$ as
$\avg{v} =  {\tfrac{1}{2}\qp{v|_{K_1} + v|_{K_2}}}$,
$\avg{\vec v} = {\tfrac{1}{2}\qp{\vec{v}|_{K_1} + \vec{v}|_{K_2}}}$ and jump operators as
$\jump{v} = {{{v}|_{K_1} \geovec n_{K_1} + {v}|_{K_2}} \geovec n_{K_2}}$,
$\jump{\vec v} = {{\vec{v}|_{K_1}}} \cdot \geovec n_{K_1} + {{\vec{v}|_{K_2}}} \cdot \geovec n_{K_2}$.
Note that on the boundary of the domain $\partial\W$ the jump and
average operators are defined as
$\avg{v} \Big\vert_{\partial\W} := v$,
$\avg{\geovec v}   \Big\vert_{\partial\W}
    :=
    \geovec v$,
$    \jump{v}
    \Big\vert_{\partial\W}
    := 
    v\geovec n$,
$    \jump{\geovec v}
    \Big\vert_{\partial\W} 
    :=
             {\geovec v}\cdot\geovec n$,

\begin{Defn}[mesh dependent norms]
  We introduce the {mesh dependent} $\leb{2}-$, $\sobh1$-and $\sobh2$-norms to be 
  \begin{gather}
    \znorm{w_h}^2 := \Norm{w_h}_{\leb{2}(\W)}^2 + \Norm{h^{3/2}\avg{\nabla w_h}}_{\leb{2}(\E\cup\partial\W)}^2 + \Norm{h^{1/2}\jump{w_h}}_{\leb{2}(\E\cup\partial\W)}^2
    \\
    \enorm{w_h}^2 := \Norm{\nabla_h w_h}_{\leb{2}(\W)}^2 + \Norm{h^{-1/2}\jump{w_h}}_{\leb{2}(\E\cup\partial\W)}^2
    \\
    \eenorm{w_h}^2 
    :=
    \Norm{\Delta_h w_h}_{\leb{2}(\W)}^2
    +
   \Norm{h^{-1/2}\jump{\nabla w_h}}_{\leb{2}(\E)}^2
    +
 \Norm{ h^{-3/2}\jump{w_h}}_{\leb{2}(\E\cup\partial\W)}^2.
  \end{gather}
  Note for $w_h \in \fes$ in view of scaling each mesh dependent norm is equaivalent to the continuous counterpart, that is $\znorm{w_h} \sim \Norm{w_h}_{\leb{2}(\W)}$ for example.
\end{Defn}
Consider the interior penalty (IP) discretisation of
(\ref{eq:bilinear-form}), to find $u_h \in \fes$ such that
\begin{equation}
  \bih{u_h}{v_h} = \ltwop{f}{v_h} \Foreach v_h \in \fes,
\end{equation}
where  
\begin{equation}
  \label{eq:IP}
  \begin{split}
    \bih{u_h}{v_h} 
    &=
    \int_\W \nabla_h u_h \cdot \nabla_h v_h \dx
    -
    \int_{\E\cup\partial\W} \jump{u_h} \cdot \avg{\nabla v_h}+\jump{v_h} \cdot \avg{\nabla u_h}\dS\\
    &\qquad +
    \int_{\E\cup\partial\W}
        \frac{\sigma_0}{h} \jump{u_h}\cdot \jump{v_h}\dS
+    \int_\E \sigma_1 h \jump{\nabla u_h}\cdot \jump{\nabla v_h}\dS
    ,
  \end{split}
\end{equation}
where $\sigma_0, \sigma_1 > 0$ represent
\emph{penalty parameters}. Note that a standard choice is to take
$\sigma_1 = 0$. The choice $\sigma_1 \neq 0$ results in a class of \emph{stabilised} dG methods \cite{Burman:2005}.
\begin{Pro}[{continuity and coercivity of $\bih{\cdot}{\cdot}$ \cite[c.f.]{ArnoldBrezziCockburnMarini:2001,ErnGuermond:2004}}]
  For $\sigma_1 \ge 0$ and $\sigma_0$ large enough and any $u_h,v_h\in\fes$ the bilinear form $\bih{\cdot}{\cdot}$ satisfies
  \begin{gather}
    \bih{u_h}{u_h} \geq C \enorm{u_h}
    \\
    \bih{u_h}{v_h} \leq C \enorm{u_h} \enorm{v_h}.
  \end{gather}
  Lax-Milgram Theorem guarantees a unique solution to the problem (\ref{eq:IP}). Also, since $u\in\sobh{2}(\W)$, the bilinear form is \emph{consistent}, hence, Strang's Lemma yields quasioptimal convergence of the method in the $\enorm{\cdot}$ norm:
  \begin{equation}
    \label{eq:quasibesth1}
    \enorm{u - u_h} \leq C \inf_{w_h\in\fes} \enorm{u - w_h}.
  \end{equation}
\end{Pro}

\textbf{Conforming reconstruction operators:} The key tool in the proof of the inf-sup condition is the notion of \emph{reconstruction operators}. It is commonplace in the a posteriori analysis of nonconforming schemes to make use of such operators. A simple, quite general methodology for the construction of
reconstruction operators is to use an averaging interpolation
operator into an $\sobh{2}$-conforming finite element
space. For example a $\cont{1}$ Hsieh--Clough--Tocher (HCT)
macro-element conforming space for $\sobh{2}$ conformity
\cite[c.f.]{BrennerGudiSung:2010,GeorgoulisHoustonVirtanen:2011}. Another option is the use of Argyris-type reconstructions \cite{BrennerGudiSung:2010}.

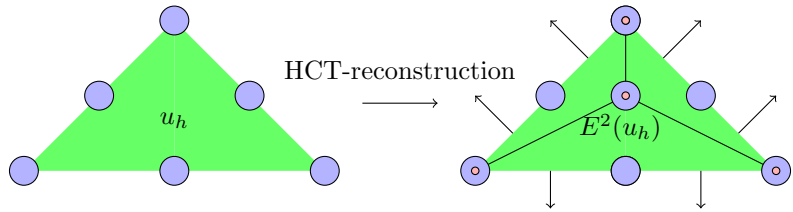
\begin{figure}[h]
  \begin{center}
    \begin{tikzpicture}      
      
      \path[fill=green!60] 
      (0,0) node (N0) {}
      -- (2,0) node (N5) {}
      -- (2,2) node (N2) {}
      -- cycle;
      
      \path[coordinate]
      (2,.67) node (K) {};
      
      \path[fill=green!60,xshift=6cm] 
      (0,0) node (NL1) {}
      -- (2,0) node (NL2) {}
      -- (2,2) node (NL0) {}
      -- cycle;
      
      \foreach \s in {0,1}{
	\path[coordinate,xshift=6cm*\s]
	(1.33,.67) node (L\s) {};
      }
      
      \path[fill=green!60]
      (4,0) node (N1) {}
      -- (2,2)
      -- (2,0)
      -- cycle
      ;
      
      \path[fill=green!60,xshift=6cm]
      (4,0) node (NR0) {}
      -- (2,2) node (NR1) {}
      -- (2,0) node (NR2) {}
      -- cycle
      ;
      
      \foreach \s in {0,1}{
	\path[coordinate,xshift=6cm*\s]
	(2.67,.67) node (R\s) {};
      }
      
      \path[coordinate]
      (3,1) node (N3) {}
      (1,1) node (N4) {}
      ;
      
      \path[coordinate,xshift=6cm]
      (1,0) node (NL3) {}
      (2,1) node (NL4) {}
      (1,1) node (NL5) {}
      ;

      \path[coordinate,xshift=6cm]
      (2,1) node (NR3) {}
      (3,0) node (NR4) {}
      (3,1) node (NR5) {}
      ;      
      
       \draw[-] (8,2) -- (8,1);
       \draw[-] (8,1) -- (6,0);
       \draw[-] (8,1) -- (10,0);
       
      \foreach \s in {0,1,...,5}{
	\node (K\s) at (N\s) [circle,inner sep=1pt,fill=blue!30,draw] {\phantom{1}};
      }

      \foreach \s in {0,1,2,5}{
	\node (KL\s) at (NL\s) [circle,inner sep=1pt,fill=blue!30,draw] {\phantom{1}};
	\node (KR\s) at (NR\s) [circle,inner sep=1pt,fill=blue!30,draw] {\phantom{1}};
      }
      \foreach \s in {0,1}{
	\node (KL\s) at (NL\s) [circle,inner sep=1pt,fill=red!30,draw] {};
	\node (KR\s) at (NR\s) [circle,inner sep=1pt,fill=red!30,draw] {};
      }

      \path[coordinate,xshift=6cm]
      (2,1) node (C1) {}
      ;      
      
      \node (K1) at (C1) [circle,inner sep=1pt,fill=blue!30,draw] {\phantom{1}};
      \node (K1) at (C1) [circle,inner sep=1pt,fill=red!30,draw] {};

      \draw[->] (7,0) -- (7,-0.5);
      \draw[->] (9,0) -- (9,-0.5);
      
      \draw[->] (6.5,0.5) -- (6,1);
      \draw[->] (7.5,1.5) -- (7,2);

      \draw[->] (9.5,0.5) -- (10,1);
      \draw[->] (8.5,1.5) -- (9,2);

      
      \node at (K) [rectangle,inner sep=1pt] {$u_h$};
      \node at (L1) [rectangle,inner sep=1pt,xshift=.6cm,yshift=-.1cm ] {$E^2(u_h)$};

      \foreach \s in {0,1}{
      }
      
      
      \draw[->] (4.5,0.9) -- (5.5,0.9)
      node[pos=0.5,above,yshift=.2cm] {HCT-reconstruction};
    \end{tikzpicture}
\end{center}
\caption{The $\poly{4}$ Hsieh--Clough--Tocher-type macro-element, as a $\sobh2(\W)$-conforming reconstruction to the quadratic Lagrange element.}
\label{fig:h2conf}
\end{figure}

\begin{Example}[$\sobh{2}(\W)$-reconstructions]
  An example of $\sobh{2}(\W)$ reconstruction operator $E^2(u_h)$ of quadratic Lagrange elements ($k=2$) is defined as follows. Let $\vec x$ be a degree of freedom of the $H^2$-conforming space $\HCT$ consisting of HCT-type macro-elements of degree $k+2$,  and let $\patch{K_{\vec x}}$ be the
  set of all elements sharing the degree of freedom $\vec x$. Then, the
  reconstruction at that specific degree of freedom is given by
  \begin{equation}
    E^2(u_h)(\vec x) 
    =
    \frac{1}{{
    \rm card}(\patch{K_{\vec x}})}\sum_{K \in \patch{K_{\vec x}}}  u_h\vert_K(\vec x).
  \end{equation}
  For the case $k=2$, the associated degrees of freedom are illustrated in Figure \ref{fig:h2conf}. Notice that the degrees of freedom of the reconstruction are a superset of those of the original finite element. This is due to the lack of existence of a conforming $\sobh2(\W)$ subspace in $\fes$ for low $k$; for instance, the existence of an $\sobh2(\W)$-conforming space requires $k\ge5$ in two dimensions (Argyris space).   Corresponding reconstructions for higher polynomial degrees have been considered in \cite{BrennerGudiSung:2010,GeorgoulisHoustonVirtanen:2011}, for instance.
\end{Example}

\begin{Lem}[{reconstruction bounds \cite[Lem 3.1]{GeorgoulisHoustonVirtanen:2011}}]
  \label{lem:reconstruction-bounds}
  The $\HCT$ reconstruction operator $E^2 : \fes\to\sobh{2}(\W)$ satisifies the following bound for all $u_h\in\fes$:
  \begin{equation}
    \enorm{E^2(u_h) - u_h}^2
    \leq
    C \qp{
       \Norm{h^{1/2}\jump{\nabla u_h}}_{\leb{2}(\E)}^2
      + 
       \Norm{h^{-1/2}\jump{u_h}}_{\leb{2}(\E)}^2
    },
  \end{equation}
with the constant $C>0$ independent of $u_h$ and of $h$.
\end{Lem}

Using this $\HCT$-reconstruction, we can construct a further $\HCT$-reconstruction admitting the same bounds, but also satisfying an \emph{adjoint orthogonality} property.
\begin{Defn}[$HCT(k+2)$-Ritz reconstruction]
  We define the Hsieh--Clough--Tocher $\sobh{2}(\W)$-conforming Ritz reconstruction operator $\R : \fes \to \HCT$ such that
  \begin{equation}
    \label{eq:hct-ritz-recon}
    \int_\W \nabla \R(u_h) \cdot \nabla \tv \dx
    =
    \int_\W \nabla_h u_h \cdot \nabla \tv\dx
    -
    \int_{\E\cup\partial\W} \jump{u_h} \cdot \nabla \tv\dS \Foreach \tv \in \HCT.
  \end{equation}
\end{Defn}

\begin{Lem}[Properties of $\R$]
  \label{lem:recon-prop}
  The $\HCT$-Ritz reconstruction is well-defined and satisfies the orthogonality condition:
  \begin{equation}
    \label{eq:orthog}
    \int_\W \qp{u_h - \R(u_h)} \Delta \tv = 0 \Foreach \tv\in\HCT.
  \end{equation}
  In addition, for $\alpha=0,1,2$, we have 
  \begin{equation}
    \label{eq:hctbd}
    \sum_{K\in\T{}}\norm{{\R(u_h) - u_h}}_{\sobh{\alpha}(K)}^2
    \leq
    C \qp{
      \Norm{h^{3/2-\alpha}\jump{\nabla u_h}}_{\leb{2}(\E)}^2
      + 
      \Norm{h^{1/2-\alpha} \jump{u_h}}_{\leb{2}(\E)}^2
    },
  \end{equation}
  for $C>0$ constants, independent of $u_h$ and of $h$.
\end{Lem}
\begin{Proof}
Fixing $u_h\in\fes$,  $\R(u_h)$ is well-defined. Indeed, setting $\tv=\R(u_h)$ in (\ref{eq:hct-ritz-recon}), along with a standard inverse estimate, we deduce 
\[
\Norm{\nabla \R(u_h)}_{\leb{2}(\W)}\le C\enorm{\nabla u_h},
\] 
for $C>0$ independent of $u_h$.
The orthogonality condition follows from integrating both sides of (\ref{eq:hct-ritz-recon}) by parts.

  To see (\ref{eq:hctbd}) we note that
  \begin{equation}
    \begin{split}
      \enorm{{\R(u_h) - u_h}}^2
      &\le
      C\bih{\R(u_h) - u_h}{\R(u_h) - u_h}
      \\
      &=
      C\bih{\R(u_h) - u_h}{E^2(u_h) - u_h} 
      \\
      &\leq
      C\enorm{\R(u_h) - u_h}
      \enorm{E^2(u_h) - u_h}.
    \end{split}
  \end{equation}
  Using the properties of $E^2(u_h)$ from Lemma \ref{lem:reconstruction-bounds} shows the claim for
  $\alpha=1$. The result for $\alpha =
    2$ follows by an inverse inequality. 
    
    For $\alpha = 0$ we use a duality argument. Take
  $z\in\sobh{2}(\W)\cap\hoz(\W)$ as the solution of the dual problem
  \begin{equation}
    -\Delta z = \R(u_h) - u_h
  \end{equation}
  then
  \begin{equation}
    \begin{split}
      \Norm{\R(u_h) - u_h}_{\leb{2}(\W)}^2
      &=
      \int_\W -\Delta z \qp{\R(u_h) - u_h}\dx
      \\
      &=
      \int_\W -\qp{\Delta z-\Delta \widetilde{z}} \qp{\R(u_h) - u_h}\dS \Foreach \widetilde{z} \in\HCT,
      \end{split}
  \end{equation}
  in view of the orthogonality property (\ref{eq:orthog}). Integrating by parts we see
  \begin{equation}
    \begin{split}
      \Norm{\R(u_h) - u_h}_{\leb{2}(\W)}^2
      &=
      \int_\W \qp{\nabla z-\nabla \widetilde{z}} \cdot \qp{\nabla \R(u_h) - \nabla_h u_h}\dx
      +
      \int_\E  \avg{\nabla z - \nabla \widetilde{z}} \cdot\jump{u_h}\dS
      \\
      &\leq
      \Norm{h^{-1}(\nabla z-\nabla \widetilde{z})}_{\leb{2}(\W)} \Norm{h(\nabla \R(u_h) - \nabla_h u_h)}_{\leb{2}(\W)} 
      \\
      &\qquad +
      \Norm{h^{-1/2}\avg{\nabla z - \nabla \widetilde{z}}}_{\leb{2}(\E)}\Norm{h^{1/2}\jump{u_h}}_{\leb{2}(\E)} .
      \end{split}
    \end{equation}
    The result follows using the approximability of the $\HCT$ space \cite{Ciarlet:1978}
    and the regularity of the dual problem. 
\end{Proof}

\begin{The}[inf--sup stability over $W(h)$]
  \label{the:inf-sup-W}
  For polynomial degree $k \geq  2$ there exists a $\gamma_h>0$, independent of $h$, such that when $\sigma_0, \sigma_1 \gg 1$
  \begin{equation}
    \label{eq:infsup2}
    \sup_{\tv \in W(h)} \frac{\bih{w_h}{\tv}}{\znorm{\tv}} \geq \gamma_h \eenorm{w_h},
  \end{equation}
  where $W(h):= \fes+\HCT$.
\end{The}

\begin{Proof}
  The proof consists of two steps. We first show there exists a $\tv\in W(h)$ such that
  \begin{equation}
    \label{eq:infsup-pf1}
    \bih{w_h}{\tv} \geq C (\min_{x\in\W}h)^2 \eenorm{w_h}^2
  \end{equation}
  and then show that 
  \begin{equation}
    \label{eq:infsup-pf2}
    \Norm{\tv}_{\leb{2}(\W)} \leq C(\max_{x\in\W}h)^2 \eenorm{w_h},
  \end{equation}
  along with the quasi-uniformity assumption on the mesh.
  
  Firstly note that, after an integration by parts, the IP method (\ref{eq:IP}) can be written as
  \begin{equation}
    \begin{split}
      \bih{w_h}{v_h}
      &=
      \int_\W - \Delta_h w_h v_h\dx
      + 
      \int_\E \jump{\nabla w_h} \avg{v_h}\dS
      - \int_{\E\cup\partial\W}
       \jump{w_h}\avg{\nabla v_h}\dS
      \\
      &\qquad +
      \int_{\E\cup\partial\W}
      \frac{\sigma_0}{h} \jump{w_h}\cdot \jump{v_h}\dS
      +      
     \int_\E \sigma_1 h \jump{\nabla w_h}\cdot \jump{\nabla v_h}\dS.
    \end{split}
  \end{equation}
Upon setting $\tv = w_h - \R(w_h) - \alpha h^2 \Delta_h w_h$, for some parameter $\alpha\in\reals$ to be chosen below, we compute 
  \begin{equation}
    \begin{split}
      \bih{w_h}{\tv}
      &=
      \alpha \Norm{h \Delta_h w_h}_{\leb{2}(\W)}^2
      +
      \sigma_1  \Norm{h^{1/2}\jump{\nabla w_h}}_{\leb{2}(\E)}^2
      +
      \sigma_0 \Norm{h^{-1/2} \jump{ w_h}}_{\leb{2}(\E\cup\partial\W)}^2
      \\
      &\qquad
      -
      \int_\W \Delta_h w_h \qp{w_h - \R(w_h)}\dx\\
      &\qquad
      +
      \int_\E \jump{\nabla w_h} \avg{w_h - \R(w_h)}\dS
      -
       \int_{\E\cup\partial\W} \jump{w_h} \cdot \avg{\nabla \qp{w_h - \R(w_h)}}\dS
      \\
      &\qquad 
      - 
      \alpha  \int_\E h^2\jump{\nabla w_h} \avg{\Delta w_h}\dS
      +
      \alpha  \int_{\E\cup\partial\W}  h^2\jump{w_h} \cdot \avg{\nabla \Delta w_h}\dS
      \\
      &\qquad 
      -\alpha 
      \sigma_1  \int_\E h^3\jump{\nabla w_h} \jump{\nabla \Delta w_h}\dS
      -
      \sigma_0 \alpha \int_{\E\cup\partial\W}   h\jump{w_h} \cdot \jump{\Delta w_h}\dS.
    \end{split}
  \end{equation}
  The orthogonality property of the $\HCT$-Ritz reconstruction (\ref{eq:hct-ritz-recon}) yields
  \begin{equation}
    \int_\W -\Delta_h w_h \qp{w_h - \R(w_h)}\dx
    =
    \int_\W \qp{\Delta \R(w_h) -\Delta_h w_h} \qp{w_h - \R(w_h)}\dx.
  \end{equation}
Repeated use of the Cauchy--Schwarz inequality, therefore, gives
  \begin{equation}
    \label{eq:pf-cs}
    \begin{split}
      \bih{w_h}{\tv}
      &\geq
      \alpha  \Norm{h\Delta_h w_h}_{\leb{2}(\W)}^2
      +
      \sigma_1  \Norm{h^{1/2}\jump{\nabla w_h}}_{\leb{2}(\E)}^2
      +
      \sigma_0  \Norm{h^{-1/2}\jump{ w_h}}_{\leb{2}(\E\cup\partial\W)}^2
      \\
      &\qquad
      -
      \Norm{h\qp{\Delta_h w_h - \Delta \R(w_h)}}_{\leb{2}(\W)}
      \Norm{h^{-1}\qp{w_h - \R(w_h)}}_{\leb{2}(\W)}
      \\
      &\qquad -
      \Norm{h^{1/2}\jump{\nabla w_h}}_{\leb{2}(\E)}
      \Norm{h^{-1/2}\avg{w_h-\R(w_h)}}_{\leb{2}(\E)}
      \\
      &\qquad - 
      \Norm{h^{-1/2}\jump{w_h}}_{\leb{2}(\E\cup\partial\W)}
      \Norm{h^{1/2}\avg{\nabla w_h-\nabla\R(w_h)}}_{\leb{2}(\E\cup\partial\W)}
      \\
      &\qquad - 
      \alpha  \Norm{h^{1/2}\jump{\nabla w_h}}_{\leb{2}(\E)}
      \Norm{h^{3/2}\avg{\Delta w_h}}_{\leb{2}(\E)}
      \\
      &\qquad 
      -
       \alpha  \Norm{h^{-1/2}\jump{w_h}}_{\leb{2}(\E\cup\partial\W)}
      \Norm{h^{5/2}\avg{\nabla \Delta w_h}}_{\leb{2}(\E\cup\partial\W)}
      \\
      &\qquad -
      \sigma_1 \alpha \Norm{h^{1/2}\jump{\nabla w_h}}_{\leb{2}(\E)}
      \Norm{h^{5/2}\jump{\nabla \Delta w_h}}_{\leb{2}(\E)}
      \\
      &\qquad 
      -
      \sigma_0 \alpha  \Norm{h^{-1/2}\jump{w_h}}_{\leb{2}(\E\cup\partial\W)}
      \Norm{h^{3/2}\jump{\Delta w_h}}_{\leb{2}(\E\cup\partial\W)}
      \\
      & =:
      \alpha \Norm{h\Delta_h w_h}_{\leb{2}(\W)}^2
      +
      \sigma_1  \Norm{h^{1/2}\jump{\nabla w_h}}_{\leb{2}(\E)}^2
      +
      \sigma_0  \Norm{h^{-1/2}\jump{w_h}}_{\leb{2}(\E\cup\partial\W)}^2
     -
      \sum_{i=1}^7 \cI_i.
    \end{split}
  \end{equation}
  We proceed to bound each of the terms $\cI_i$ individually. Note
  that in view of scaling and inverse inequalities we have for any
  $w_h\in\fes$:
  \begin{gather}
    \label{eq:useful-1}
    \Norm{\avg{w_h}}_{\leb{2}(e)}
    \leq C_1  \Norm{h^{-1/2}w_h}_{\leb{2}(\bar{K}_1\cup\bar{K}_2)}
    \\
    \label{eq:useful-2}
    \Norm{\avg{\nabla w_h}}_{\leb{2}(e)}
    \leq C_2 \Norm{ h^{-3/2} w_h}_{\leb{2}(\bar{K}_1\cup\bar{K}_2)}
  \end{gather}
  for any edge/face $e:=\bar{K}_1\cap \bar{K}_2\in \E$, and elements
  $K_1, K_2\in \T{}$, with $C_1,C_2$ depending only on the mesh-regularity
  and shape-regularity constants.
  
  For $\cI_1$, in view of Lemma \ref{lem:recon-prop}, we have
  \begin{equation}
    \label{eq:pf1}
    \begin{split}
      \cI_1 
      &\leq
      {C_3} \qp{\Norm{h^{1/2} \jump{\nabla w_h}}_{\leb{2}(\E)}^2 +  \Norm{h^{-1/2} \jump{w_h}}_{\leb{2}(\E)}^2},
    \end{split}
  \end{equation}
with constant $C_3>0$ being the maximum of all constants in (\ref{eq:hctbd}) for all $\alpha$.
  
  For $\cI_2$, (\ref{eq:useful-1}) and Lemma  \ref{lem:recon-prop} yield
  \begin{equation}
    \begin{split}
      \cI_2
      &\leq 
      C_1 {C_3^{1/2}} 
      \Norm{h^{1/2}\jump{\nabla w_h}}_{\leb{2}(\E)}
      \qp{\Norm{h^{1/2} \jump{\nabla w_h}}_{\leb{2}(\E)}^2 +  \Norm{h^{-1/2} \jump{w_h}}_{\leb{2}(\E)}^2}^{1/2}
      \\
      &\leq
      C_1^2 C_3 
      \Norm{h^{1/2}\jump{\nabla w_h}}_{\leb{2}(\E)}^2
      +
      \frac{C_1^2 C_3 }{2}
      \Norm{h^{-1/2}\jump{w_h}}_{\leb{2}(\E)}^2.
    \end{split}
  \end{equation}
 
  For  $\cI_3$, (\ref{eq:useful-2}) and Lemma  \ref{lem:recon-prop} yield
  \begin{equation}
    \begin{split}
      \cI_3
      &\leq
      C_2C_3^{1/2}  \Norm{h^{-1/2}\jump{w_h}}_{\leb{2}(\E\cup\partial\W)}
      \qp{\Norm{h^{1/2} \jump{\nabla w_h}}_{\leb{2}(\E)}^2 +  \Norm{h^{-1/2} \jump{w_h}}_{\leb{2}(\E)}^2}^{1/2}
      \\
      &\leq
      C_2^2C_3
      \Norm{h^{1/2}\jump{\nabla w_h}}_{\leb{2}(\E\cup\partial\W)}
      +
      \frac{C_2^2C_3 }{2} \Norm{h^{-1/2}\jump{w_h}}_{\leb{2}(\E)}.
    \end{split}
  \end{equation}
  
  For $\cI_4$, we have
  \begin{equation}
    \begin{split}
      \cI_4
      \leq 
      C_1 \alpha  \Norm{h^{1/2}\jump{\nabla w_h}}_{\leb{2}(\E)}\Norm{h\Delta w_h}_{\leb{2}(\W)}
      \leq
      \epsilon_4 \alpha \Norm{h\Delta w_h}_{\leb{2}(\W)}^2
      +
      \frac{C_1^2\alpha}{4\epsilon_4} \Norm{h^{1/2}\jump{\nabla w_h}}_{\leb{2}(\E)}^2,
    \end{split}
  \end{equation}
  for any $\epsilon_4>0$, while for $\cI_5$, we get
  \begin{equation}
    \begin{split}
      \cI_5
      \leq
      C_2  \alpha  \Norm{h^{-1/2}\jump{w_h}}_{\leb{2}(\E\cup\partial\W)} \Norm{h\Delta w_h}_{\leb{2}(\W)}
      \leq
      \epsilon_5 \alpha  \Norm{h\Delta w_h}_{\leb{2}(\W)}^2
      +
      \frac{C_2^2\alpha }{4\epsilon_5}  \Norm{h^{-1/2}\jump{w_h}}_{\leb{2}(\E\cup\partial\W)}^2.
    \end{split}
  \end{equation}
  for any $\epsilon_5>0$; similarly for $\cI_6$ and for any $\epsilon_6>0$, we have
  \begin{equation}
    \label{eq:pf6}
    \begin{split}
      \cI_6
      \leq
      C_2 \sigma_1 \alpha \Norm{ h^{1/2}\jump{\nabla w_h}}_{\leb{2}(\E)} \Norm{h\Delta w_h}_{\leb{2}(\W)}
      \leq
      \epsilon_6 \alpha   \Norm{h\Delta w_h}_{\leb{2}(\W)}^2 
      +
      \frac{C_2^2 \sigma_1^2 \alpha}{4\epsilon_6}  \Norm{h^{1/2}\jump{\nabla w_h}}_{\leb{2}(\E)}^2.
    \end{split}
  \end{equation}
  Finally, the last term $\cI_7$ can be bounded as follows:
  \begin{equation}
    \label{eq:pf7}
      \cI_7
      \leq
      C_1 \sigma_0 \alpha \Norm{ h^{-1/2}\jump{w_h}}_{\leb{2}(\E\cup\partial\W)} \Norm{h\Delta_h w_h}_{\leb{2}(\W)}
      \leq
      \epsilon_7 \alpha  \Norm{h\Delta_h w_h}_{\leb{2}(\W)}^2
      +
      \frac{C_1^2 \sigma_0^2 \alpha}{4\epsilon_7} \Norm{h^{-1/2} \jump{w_h}}_{\leb{2}(\E\cup\partial\W)}^2,
  \end{equation}
  for any $\epsilon_7>0$.
  
  Collecting the results (\ref{eq:pf1})--(\ref{eq:pf7}) and
  substituting this into (\ref{eq:pf-cs}) we deduce
  \begin{equation}
    \begin{split}
      \bih{w_h}{\tv}
      &\geq
       \Norm{h\Delta_h w_h}_{\leb{2}(\W)}^2
     \alpha  \Big(1 - \epsilon_4 -  \epsilon_5 -  \epsilon_6 - \epsilon_7\Big)
      \\
      &\qquad +
      \Norm{h^{1/2}\jump{\nabla w_h}}_{\leb{2}(\E)}^2
      \Big(
      \sigma_1
      -
      C_3
      -
      C_1^2C_3
      -
      \frac{C_2^2C_3 ^2}{2}
      -
      \frac{C_1^2 \alpha}{4\epsilon_4}
      -
      \frac{C_2^2 \sigma_1^2 \alpha}{4\epsilon_6}
      \Big)
      \\
      &\qquad 
      +
      \Norm{h^{-1/2}\jump{w_h}}_{\leb{2}(\E\cup\W)}^2
      \Big(
      \sigma_0
      -
      C_3
      -
      \frac{C_1^2C_3}{2}
      -
      C_2^2C_3
      -
      \frac{C_2^2 \alpha }{4\epsilon_5}
      -
      \frac{C_1^2 \sigma_0^2 \alpha}{4\epsilon_7}
      \Big).
    \end{split}
  \end{equation}
  To arrive to (\ref{eq:infsup-pf1}), we can choose $\epsilon_4 = \epsilon_5
  =\epsilon_6 = \epsilon_7 = \frac{1}{5}$, $\alpha =(\max\qp{\sigma_0^2,\sigma_1^2})^{-1}$ and $\sigma_0$ and
  $\sigma_1$ large enough.

For (\ref{eq:infsup-pf2}), we use Lemma \ref{lem:recon-prop} to see that
  \begin{equation}
    \Norm{\tv}_{\leb{2}(\W)}
    \leq
    \Norm{w_h - \R(w_h)}_{\leb{2}(\W)} + \Norm{\alpha h^2 \Delta_hw_h}_{\leb{2}(\W)}
    \leq
    C \eenorm{h^2 w_h},
  \end{equation}
  yielding (\ref{eq:infsup2}).
\end{Proof}

\begin{Lem}[Stability of the Ritz projection]
  \label{lem:stabilityofRitz}
  Let $R$ demote the $\bih{\cdot}{\cdot}$ orthogonal projector into $\fes$, then for $\tw \in W(h)$ we have that
  \begin{equation}
    \Norm{R \tw}_{\leb{2}(\W)}
    \leq
    C
    \qp{
      \Norm{h\nabla \tw}_{\leb{2}(\W)}
      +
      \Norm{\tw}_{\leb{2}(\W)}
    }
    \leq
    C \Norm{\tw}_{\leb{2}(\W)}.
  \end{equation}
\end{Lem}
\begin{Proof}
  Let $\tg\in W(h)$ be the solution to the discrete dual problem such
  that
  \begin{equation}
    \begin{split}
      \bih{\tv}{\tg} &= \ltwop{R \tw}{\tv} \Foreach \tv \in W(h).
    \end{split}
  \end{equation}
  Then we have
  \begin{equation}\label{eq:duality_one}
    \begin{split}
      \Norm{R \tw}_{\leb{2}(\W)}^2
      &=
      \ltwop{R\tw - \tw}{R\tw} + \ltwop{\tw}{R\tw}
      \\
      &=
      \bih{R \tw - \tw}{\tg} + \ltwop{\tw}{R\tw}.
    \end{split}
  \end{equation}
  Let $\Pi: \sobh{1}(\W)\to \fes\cap \sobh{1}_0(\W)$ a suitable projection with optimal approximation properties. Then
  \begin{equation}\label{eq:duality_two}
    \begin{split}
      \Norm{R \tw}_{\leb{2}(\W)}^2
      &=
      \bih{R \tw - \tw}{\tg - \Pi \tg} + \ltwop{\tw}{R\tw}
      \\
      &\le
      \enorm{h (R\tw-\tw)}
      \enorm{h^{-1}(\tg-\Pi \tg)}
      +
      \Norm{\tw}_{\leb{2}(\W)} \Norm{R\tw}_{\leb{2}(\W)},
    \end{split}
  \end{equation}
  through the continuity of $\bih{\cdot}{\cdot}$.
  From the optimal approximation properties of the projection/interpolant $\Pi$, we have
 \begin{equation}
 \label{eq:duality_three}
  \begin{split}
  &\enorm{h^{-1}\qp{ \tg - \Pi \tg}}^2
  \le  \tilde C\eenorm{\tg}^2,
  \end{split}
 \end{equation}
 and, using the discrete regularity of $\tg$ induced by the inf-sup
 condition in Theorem \ref{the:inf-sup-W}
 \begin{equation}
   \gamma_h\eenorm{\tg}
   \leq
   \sup_{\tv\in W(h)} \frac{\bih{\tg}{\tv}}{\znorm{\tv}}
   =
   \sup_{\tv\in W(h)} \frac{\ltwop{R \tw}{\tv}}{\znorm{\tv}}
   \leq
   C \Norm{R\tw}_{\leb{2}(\W)}.
 \end{equation}
 Hence we see that
 \begin{equation}
   \begin{split}
     \Norm{R \tw}_{\leb{2}(\W)}^2
     &\leq
     C\qp{
       \enorm{h (R\tw-\tw)}
     \Norm{R\tw}_{\leb{2}(\W)}
     +
     \Norm{\tw}_{\leb{2}(\W)} \Norm{R\tw}_{\leb{2}(\W)}
     }
     \\
     &\leq
     \enorm{h\tw}
     \Norm{R\tw}_{\leb{2}(\W)}
     +
     \Norm{\tw}_{\leb{2}(\W)} \Norm{R\tw}_{\leb{2}(\W)},
   \end{split}
 \end{equation}
 in view of the quasi-best approximation in $\enorm{\cdot}$ from (\ref{eq:quasibesth1}). The conclusion follows from standard inverse inequalities.
\end{Proof}

\begin{The}[inf--sup stability over $\fes$]
  \label{the:inf-sup-V}
  For polynomial degree $k \geq  2$ there exists a $\gamma_h>0$, independent of $h$, such that when $\sigma_0, \sigma_1 \gg 1$
  \begin{equation}
    \label{eq:infsup1}
    \sup_{v_h\in\fes} \frac{\bih{w_h}{v_h}}{\eenorm{v_h}} \ge \gamma_h \znorm{w_h}
  \end{equation}
\end{The}

\begin{Proof}
To show (\ref{eq:infsup1}) we fix $w_h$ and let $\Phi\in\fes$ be the solution of the dual problem
  \begin{equation}
    \bih{\Psi}{\Phi} = \int_\W w_h \Psi \dx \Foreach \Psi\in\fes.
  \end{equation}
  Following the same arguments as in the proof of Theorem \ref{the:inf-sup-W}, it is clear that there exists a $C>0$ such that
  \begin{equation}
    Ch^2 \eenorm{\Phi}^2 \leq \bih{\Phi}{\tv},
  \end{equation}
  where $\tv := \Phi - \R(\Phi) - \alpha h^2 \Delta_h \Phi$. Now it is clear that
  \begin{equation}
    \Norm{\tv}_{\leb{2}(\W)} \eenorm{\Phi} \leq Ch^2 \eenorm{\Phi}^2 \leq C \bih{\Phi}{\tv},
  \end{equation}
  and hence in view of Lemma \ref{lem:stabilityofRitz}  we have, with $R$ denoting the $\cA_h$ orthogonal projector into $\fes$, that
  \begin{equation}
    \Norm{R \tv}_{\leb{2}(\W)}
    \leq
    C \qp{\Norm{h \nabla \tv}_{\leb{2}(\W)} + \Norm{\tv}_{\leb{2}(\W)}}
    \leq
    C \Norm{\tv}_{\leb{2}(\W)},
  \end{equation}
  through inverse inequalities. Hence arguing as above
  \begin{equation}
    \Norm{R \tv}_{\leb{2}(\W)}
    \eenorm{\Phi}
    \leq
    C\Norm{\tv}_{\leb{2}(\W)}
    \eenorm{\Phi} 
    \leq
    Ch^2\eenorm{\Phi}^2
    \leq
    C\bih{\Phi}{\tv}
    =
    C\bih{\Phi}{R \tv},
  \end{equation}
  concluding the proof.
\end{Proof}

\begin{Cor}[Convergence]
  \label{cor:convergence}
  Let $u$ solve (\ref{eq:bilinear-form}) and $u_h \in \fes$ be the interior penalty
  approximation from (\ref{eq:IP}), then
  \begin{equation}
    \znorm{u - u_h}
    \leq
    \qp{1+\frac{C_B}{\gamma_h}}\inf_{w_h\in\fes}
    \znorm{u - w_h}
    +
    \frac 1 {\gamma_h}
    \sup_{v_h\in\fes}
    \frac{\bih{u_h - u}{v_h}}{\eenorm{v_h}}    
    .
  \end{equation}
  If $u\in\sobh{k+1}(\W)$ the following a priori bound holds:
  \begin{equation}
    \label{eq:apriori}
    \znorm{u-u_h}
    +
    \enorm{h\qp{u-u_h}}
    +
    \eenorm{h^2\qp{u - u_h}} \leq C h^{k+1} \norm{u}_{k+1}.
  \end{equation}
\end{Cor}
\begin{Proof}
  Using the inf-sup condition from Theorem \ref{the:inf-sup-V} we see
  \begin{equation}
    \begin{split}
      \gamma_h \znorm{w_h - u_h}
      &\leq
      \sup_{v_h\in\fes} \frac{\bih{w_h - u_h}{v_h}}{\eenorm{v_h}}
      \\
      &\leq
      \sup_{v_h\in\fes} \frac{\bih{w_h - u}{v_h}}{\eenorm{v_h}}
      +
      \sup_{v_h\in\fes} \frac{\bih{u - u_h}{v_h}}{\eenorm{v_h}}.
    \end{split}
  \end{equation}
  Now using the natural continuity bound
  \begin{equation}
    \bih{u-w_h}{v_h} \leq C_B \znorm{u-w_h} \eenorm{v_h}
  \end{equation}
  we see
  \begin{equation}
    \znorm{w_h - u_h}
    \leq
    \frac{C_B}{\gamma_h} \znorm{u-w_h}
    +
    \frac{1}{\gamma_h}
    \sup_{v_h\in\fes} \frac{\bih{u - u_h}{v_h}}{\eenorm{v_h}}.
  \end{equation}
  Hence, in view of the triangle inequality
  \begin{equation}
    \znorm{u - u_h}
    \leq
    \qp{1 + \frac{C_B}{\gamma_h}} \znorm{u-w_h}
    +
    \frac{1}{\gamma_h}
    \sup_{v_h\in\fes} \frac{\bih{u - u_h}{v_h}}{\eenorm{v_h}}.
  \end{equation}
  The conclusion follows since $w_h$ was arbitrary.
\end{Proof}

\section{Applications to problems with Rough Data}
\label{sec:roughdata}
In this section we examine the problem
\begin{equation}
  \label{eq:roughprob}
  -\Delta u = f := \partial_x \delta_{\bar x},
\end{equation}
where $\delta_{\bar x}$ denotes the Dirac distribution at a point
$\bar x \in \W$. In this case we have $f\in\sobh{-2}(\W) \not
\ \ \sobh{-1}(\W)$ and hence the solution $u\in\leb{2}(\W) \not
\ \ \sobh{1}(\W)$. This means it cannot be characterised through a weak
formulation, rather an \emph{ultra weak formulation}, where we seek
$u\in\leb{2}(\W)$ such that
\begin{equation}
  \label{eq:contprob}
  \int_\W -u \Delta v \d x= \duality{f}{v}_{\sobh{-2}(\W)\times \sobhz2(\W)} \Foreach v \in \sobhz{2}(\W),
\end{equation}
and the right hand side of (\ref{eq:contprob}) is understood as a \emph{duality pairing}. 
In this setting standard methods pertaining to the analysis of
Galerkin methods may not apply, for example the Aubin--Nitsche duality
arguement. However, if we assume that $\bar x$ does not lie on the
skeleton of the triangulation the stabilised IP method is still well
defined and the inf-sup condition still holds. We define our approximation as seeking $u_h\in\fes$ such that
\begin{equation}
  \label{eq:disprob}
  \bih{u_h}{v_h} =  \duality{f}{E^2(v_h)}_{\sobh{-2}(\W)\times \sobhz2(\W)} = E^2(v_h)(\bar x) \Foreach v_h\in\fes.
\end{equation}
Since the inf-sup condition given in Theorem \ref{the:inf-sup-V} is a condition only on the operator itself the first statement in Corollary \ref{cor:convergence} holds true. The only uncertainty with the bound is the behaviour of the inconsistency term. The control of this term is the main motivation in the nonstandard definition of the right hand side of (\ref{eq:disprob}).

\begin{The}[quasi-optimal error control for problems with rough data]
  \label{the:convergence}
  Let $u\in\leb{2}(\W)$ solve (\ref{eq:contprob}) and $u_h\in\fes$ be the approximation defined through (\ref{eq:disprob}), then 
  \begin{equation}
    \Norm{u - u_h}_{\leb{2}(\W)} \leq C\inf_{w_h \in \fes} \qp{\Norm{u - w_h}_{\leb{2}(\W)} + \Norm{h^{3/2} \jump{\nabla w_h}}_{\leb{2}(\E)} + \Norm{h^{1/2} \jump {w_h}}_{\leb{2}(\E)}}.
  \end{equation}
\end{The}
\begin{Proof}
  The proof of this fact takes some inspiration from that of \cite{Gudi:2010} where inconsistency terms arise from the fact that the solution of an elliptic problem may only lie in $\sobh 1(\W)$ for which the operator $\bih{u}{v_h}$ may not be well defined. Here we have even more difficulty since the solution $u\in\leb{2}(\W)\not \ \sobh{1}(\W)$.

  Using the inf-sup condition from Theorem \ref{the:inf-sup-V} we have
  \begin{equation}
    \begin{split}
      \gamma_h \znorm{w_h - u_h}
      &\leq
      \sup_{v_h\in\fes} \frac{\bih{w_h - u_h}{v_h}}{\eenorm{v_h}}.
    \end{split}
  \end{equation}
  Now by adding and subtracting appropriate terms and using (\ref{eq:contprob}) and (\ref{eq:disprob}) we see
  \begin{equation}
    \label{pf0}
    \begin{split}
      \bih{w_h - u_h}{v_h}
      &=
      -\int_\W u \Delta \R(v_h) \d x
      -
      \bih{w_h}{\R(v_h)}
      +
      \duality{f}{E^2(v_h) - \R(v_h)}
      \\
      &\qquad -
      \bih{w_h}{E^2(v_h) - \R(v_h)}
      +
      \bih{w_h}{E^2(v_h) - v_h}
      \\
      &=
      -\int_\W u \Delta \R(v_h) \d x
      -
      \bih{w_h}{\R (v_h)}
      +
      \duality{f}{E^2(v_h) - \R(v_h)}
      \\
      &\qquad -
      \bih{w_h}{E^2(v_h) - \R(v_h)}
      +
      \bih{w_h - \R(w_h)}{E^2(v_h) - v_h},
    \end{split}
  \end{equation}
  by the orthogonality properties of $\R(w_h)$ given in (\ref{eq:hct-ritz-recon}). Now we may use that
  \begin{equation}
    \label{pf1}
    \begin{split}
      -\int_\W u \Delta \R(v_h) \d x
      -
      \bih{w_h}{\R (v_h)}
      &\leq
      C \Norm{u - w_h}_{\leb{2}(\W)}
      \Norm{\Delta \R(v_h)}_{\leb{2}(\W)}
      \\
      &\leq
      C \Norm{u - w_h}_{\leb{2}(\W)}
      \eenorm{v_h}
      ,
    \end{split}
  \end{equation}
  through the stability of $\R(v_h)$. In addition by the definition of the solution to the PDE
  \begin{equation}
    \label{pf2}
    \begin{split}
      \duality{f}{E^2(v_h) - \R(v_h)}
      -
      \bih{w_h}{E^2(v_h) - \R(v_h)}
      &=
      -
      \int_\W u \Delta \qp{E^2(v_h) - \R(v_h)} \d x
      \\
      &\qquad \qquad -
      \bih{w_h}{E^2(v_h) - \R(v_h)}
      \\
      &\leq
      C \Norm{u - w_h}_{\leb{2}(\W)}
      \Norm{\Delta \qp{E^2(v_h) - \R(v_h)}}_{\leb{2}(\W)}
      \\
      &\leq
      C \Norm{u - w_h}_{\leb{2}(\W)}
      \eenorm{v_h}.
    \end{split}
  \end{equation}
  Finally, using the approximation properties of $\R(\cdot)$ and $E^2(\cdot)$ we see
  \begin{equation}
    \label{pf3}
    \begin{split}
      \bih{w_h - \R(w_h)}{E^2(v_h) - v_h}
      &\leq
      C \enorm{w_h - \R(w_h)} \enorm{E^2(v_h) - v_h}
      \\
      &\leq
      C \qp{\Norm{h^{3/2} \jump{\nabla w_h}}_{\leb{2}(\E)} + \Norm{h^{1/2} \jump {w_h}}_{\leb{2}(\E)}} \eenorm{v_h}.
    \end{split}
  \end{equation}
  Substituting (\ref{pf1}), (\ref{pf2}) and (\ref{pf3}) into (\ref{pf0}) we have
  \begin{equation}
    \bih{w_h - u_h}{v_h}
    \leq
    C \qp{\Norm{u-w_h}_{\leb{2}(\W)} + \Norm{h^{3/2} \jump{\nabla w_h}}_{\leb{2}(\E)} + \Norm{h^{1/2} \jump {w_h}}_{\leb{2}(\E)}} \eenorm{v_h},
  \end{equation}
  and hence
  \begin{equation}
    \begin{split}
      \Norm{u - u_h}_{\leb{2}(\W)}
      &\leq
      \Norm{u - w_h}_{\leb{2}(\W)}
      +
      \Norm{u_h - w_h}_{\leb{2}(\W)}
      \\
      &\leq
      \Norm{u - w_h}_{\leb{2}(\W)}
      +
      \znorm{u_h - w_h}
      \\
      &\leq
      C \qp{\Norm{u-w_h}_{\leb{2}(\W)} + \Norm{h^{3/2} \jump{\nabla w_h}}_{\leb{2}(\E)} + \Norm{h^{1/2} \jump {w_h}}_{\leb{2}(\E)}},
    \end{split}
  \end{equation}
  as required.
\end{Proof}

\subsection{A numerical experiment}

\begin{figure}[!ht]
  \caption[]
          {\label{fig:sol}
            In this experiment we test the $\leb{2}$ convergence of
            the interior penality method to approximate the
            distributional solution (\ref{eq:dist-sol}).
         }
         \begin{center}
           \subfigure[{\label{fig:b1}
               The IP approximation.
           }]{
             \includegraphics[scale=\figscale,width=0.47\figwidth]{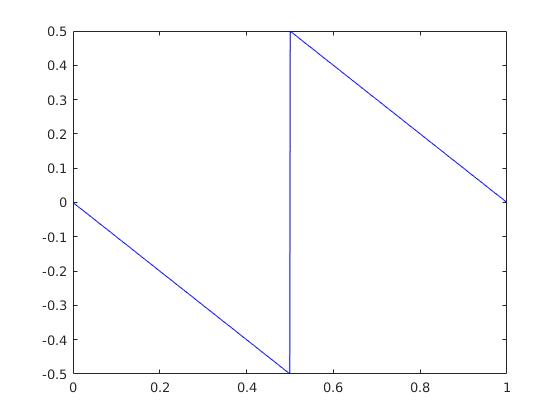}
           }
           \subfigure[{\label{fig:b1}
               The experimental order of convergence.
           }]{
             \includegraphics[scale=\figscale,width=0.47\figwidth]{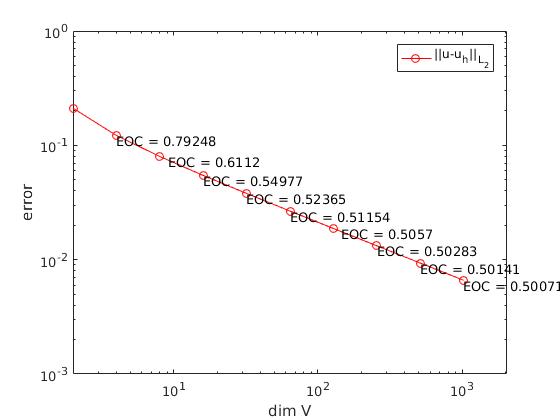}
           }
         \end{center}
 \end{figure}

If $d=1$ we can even characterise a distributional solution, indeed
for $\W = [0,1]$ we have that
\begin{equation}
  \label{eq:dist-sol}
  u(x) =
  \begin{cases}
    -x \text{ when } x < \bar x
    \\
    x \text{ when } x > \bar x,
  \end{cases}
\end{equation}
solves (\ref{eq:roughprob}). Using Theorem \ref{the:convergence} we are
able to show this approximation satisfies the a priori bound
\begin{equation}
  \Norm{u-u_h}_{\leb{2}(\W)} \leq C \sqrt{h},
\end{equation}
as verified in Figure \ref{fig:sol}.

\end{document}